\documentclass{article}

\usepackage{amssymb}

\newcommand{\proof}{{\bf Proof:  }}

\newcommand{\example}{{\bf Example:  }}

\newcommand{\hb}{\newline\hspace*{\fill}$\Box$}

\newtheorem{theorem}{Theorem}[section]
\newtheorem{lemma}[theorem]{Lemma}
\newtheorem{definition}[theorem]{Definition}
\newtheorem{proposition}[theorem]{Proposition}
\newtheorem{corollary}[theorem]{Corollary}

\begin{document}

\parindent0pt

\title{\bf Localization in quiver moduli}

\author{Markus Reineke\\ Mathematisches Institut\\ Universit\"at M\"unster\\ D - 48149 M\"unster, Germany\\
e-mail: reinekem@math.uni-muenster.de}

\date{}

\maketitle

\begin{abstract} The fixed point set under a natural torus action on
projectivized moduli spaces of simple representations of quivers is described. As an application, the Euler
characteristic of these moduli is computed.
\end{abstract}

\section{Introduction}\label{introduction}

It is a common theme in moduli theory to consider fixed point sets under torus actions, see for example
\cite[Section 3]{Kl} in the context of moduli of vector bundles.
As a general principle, the torus fixed points in moduli spaces of objects like representations, sheaves, etc.,
correspond to certain graded, or equivariant, such objects.\\[2ex]
In the ideal case of a smooth projective moduli space and a torus acting with finitely many fixed points, one can
apply the Bialynicki-Birula theorem \cite{BB} to conclude rationality, to construct cell
decompositions and to compute the Betti numbers.\\[2ex]
In general, the moduli spaces in question are seldomly projective, and the only known compactifications
are usually highly singular.
For example, this may be the
case when the moduli are constructed as Mumford quotients \cite{M} of stable points; then a natural 
compactification is often provided by considering the quotient of semistable points. In this case, the localization
to torus fixed points can still provide cohomological information on the moduli space (although the method is not strong
enough to provide cell decompositions). For example, for equivariantly formal (in the sense of \cite{GKM}) moduli,
the equivariant cohomology can be computed.\\[2ex]
The main theme of the present paper is to apply the localization principle to projectivizations
of quiver moduli parametrizing simple
representations of quivers, to be introduced and discussed in section \ref{projectivized} 
(unfortunately, there seems to be no straightforward generalization to the more general
quiver moduli
parametrizing stable representations). 
Based on an analysis of the fixed points in section \ref{fps}, the first main result,
Theorem \ref{fmr}, gives a description
of a set of torus fixed points as
a disjoint union of moduli spaces of the same type, that is, of projectived moduli of simple representations for 
other quivers. In fact, the quivers appearing in this description are `almost universal abelian coverings' of the
original quiver. Therefore, in the course of the proof, some principles of covering theory
of finite dimensional algebras have to be studied in the context of simple quiver representations
in section \ref{construction}.\\[2ex]
In \cite{RCRP}, it is proven that the number of rational points of models of quiver moduli over
finite fields behaves polynomially in the cardinality $q$ of the field. Starting from computer experiments,
a formula for the value of the first derivative of these counting polynomials at $q=1$ was conjectured.
Using some arithmetic geometry, this number can be reinterpreted as the Euler characteristic
in cohomology with compact support of the projectivized quiver moduli.\\[1ex]
In this form, the conjecture is proven in section \ref{ind}. Namely, Theorem \ref{euler} identifies
the Euler characteristic in question as the number of cyclic equivalence classes of primitive cycles
in the given quiver of a certain weight.
In fact, the localization techniques of this paper were developed for the purpose of computing this number.\\[2ex]
The description of the fixed point set provides an interesting interaction {between} algebraic geometry and representation
theory of quivers: the detection of torus fixed points naturally leads to the study of covering techniques for simple
representations. After such methods are developed, they can be successfully applied to the geometry of quiver moduli.
The description also motivates consideration of a very special class of combinatorially defined simple representations of
quivers, called string representations (Definition \ref{string}),
which are in some sense `responsible' for aspects of the global topology of the moduli.\\[3ex]
{\bf Acknowledgments:} The author would like to thank M.~Brion and W.~{Soer\-gel} for helpful remarks concerning the
localization principle, and V.~Ginzburg, S.~Kumar, L.~Le Bruyn and C.~M.~Ringel for interesting discussions.

\section{Projectivized quiver moduli}\label{projectivized}

Let $Q$ be a quiver with finite set of vertices $Q_0$ and finite set of arrows $Q_1$,
which will be written as $\alpha:i\rightarrow j$ for $i,j\in Q_0$. In the
following, we will frequently consider the abelian groups
$${\bf Z}Q_0=\bigoplus_{i\in Q_0}{\bf Z}i,\;\;\; {\bf Z}Q_1=\bigoplus_{\alpha\in Q_1}{\bf Z}\alpha,$$
and the subsemigroups ${\bf N}Q_0$ (called the lattice of dimension vectors for $Q$), resp.~${\bf N}Q_1$.
On these groups, we have natural
${\bf Z}$-gradings by setting the degree of the natural basis elements $i$ (resp.~$\alpha$) to be $1$. The degree of a
homogeneous element $d$ in either of the two groups is denoted by $|d|$.\\[1ex]
A cycle $\omega=(\alpha_1,\ldots,\alpha_s)$ at $i\in Q_0$ in $Q$ is a (possibly empty) sequence of arrows
$$i=i_0\stackrel{\alpha_1}{\rightarrow}i_1\stackrel{\alpha_2}{\rightarrow}\ldots\stackrel{\alpha_s}{\rightarrow}
i_{s}=i$$ in $Q$. To $\omega$, we associate its weight $|\omega|\in{\bf N}Q_1$ as follows: for each $\alpha\in Q_1$,
the component $|\omega|_\alpha$ equals the number of occurences of the arrow $\alpha$ in the sequence of
arrows $\alpha_1,\ldots,\alpha_s$. We associate to the cycle $\omega$ its dimension vector $\dim \omega\in{\bf N}Q_0$
as follows:
for each $i\in Q_0$, the component $(\dim\omega)_i$ equals the number of indices $l$ such that $i_l=i$.
In other words, $\omega$ has dimension vector $\sum_{i\in Q_0}d_ii$ if it passes
through each vertex $i$ precisely $d_i$-times.\\[1ex]
We have the notion of cyclic equivalence of cycles, which is the
equivalence relation generated by
$$(\alpha_1,\ldots,\alpha_s)\sim(\alpha_2,\ldots,\alpha_s,\alpha_1).$$
We denote by $C_d(Q)$ the set of cyclic equivalence classes of cycles in $Q$ of dimension vector $d$.
The cycle $\omega$ is called primitive if it is not cyclically equivalent to a non-trivial power (that is,
iterated concatenation) $(\omega')^n$ of another cycle $\omega'$. The subset of $C_d(Q)$ corresponding to classes of
primitive cycles is denoted by $C_d^{prim}(Q)$.\\[1ex]
We introduce a non-symmetric bilinear form, the Euler form, on ${\bf Z}Q_0$ by
$$\langle d,e\rangle:=\sum_{i\in Q_0}d_ie_i-\sum_{\alpha:i\rightarrow j}d_ie_j.$$
Let $$d=\sum_{i\in Q_0}d_ii\in{\bf N}Q_0$$ be a dimension vector, and fix complex vector spaces $V_i$ of
dimension $d_i$ for all $i\in Q_0$. The affine space $$R_d(Q):=\bigoplus_{\alpha:i\rightarrow j}{\rm Hom}(V_i,V_j)$$
parametrizes complex representations
$$X=((V_i)_{i\in Q_0},(X_\alpha:V_i\rightarrow V_j)_{(\alpha:i\rightarrow j)\in Q_1})$$
of $Q$ of dimension vector $d$.
The reductive linear algebraic group $$G_d:=\prod_{i\in Q_0}{\rm GL}(V_i)$$ acts on $R_d(Q)$ by base change
$$(g_i)_i\cdot(X_\alpha)_\alpha:=(g_jX_\alpha g_i^{-1})_{\alpha:i\rightarrow j}$$ in such a way that the orbits $G_dX$
correspond to the isomorphism classes $[X]$ of representations. The scalar matrices, embedded diagonally into $G_d$, act
trivially on $R_d(Q)$. Therefore, the action of $G_d$ factors through the quotient $$PG_d:=G_d/{\bf C}^*.$$

By \cite{A}, an orbit $G_dX$ is closed in $R_d(Q)$ if and only $X$ is a semisimple representation of $Q$.
Each cycle $\omega=(\alpha_1,\ldots,\alpha_s)$ in $Q$ induces a function $t_\omega$ on $R_d(Q)$
(the trace along the oriented cycle $\omega$) by
$$t_\omega:(X_\alpha)_\alpha\mapsto {\rm tr}(X_{\alpha_s}\ldots X_{\alpha_1}).$$
For each vertex $i\in Q_0$, we have the simple representation $S_i$ of dimension vector $i\in{\bf N}Q_0$.
The point $0\in R_d(Q)$ corresponds to the semisimple representation
$\bigoplus_{i\in I}S_i^{d_i}$.\\[1ex]
The following properties of a representation $X\in R_d(Q)$ are equivalent (see \cite{LBN}):

\begin{enumerate} 
\item the orbit closure $\overline{G_dX}$ contains the point $0\in R_d(Q)$, 
\item All functions $t_\omega$ vanish at the point $X$,
\item all Jordan-H\"older factors of $X$ are isomorphic to some $S_i$.
\end{enumerate}

Denote by $N_d(Q)$ the closed subvariety of $R_d(Q)$ of all points $X$ satisfying one of the above three equivalent
conditions, called the nullcone in $R_d(Q)$.\\[1ex]
In order to form the quotient of $R_d(Q)$ by the action of $G_d$, we have to describe the ring of
$G_d$-invariant functions on $R_d(Q)$. This is done by the following theorem, due to L.~Le Bruyn and C.~Procesi:

\begin{theorem}\cite{LBP1}\label{lbp} The invariant ring ${\bf C}[R_d(Q)]^{G_d}$
is generated by the traces $t_\omega$ along cycles $\omega$ in $Q$.
\end{theorem}

Therefore, one can coordinatize the quotient
$$M_d^{ssimp}(Q):=R_d(Q)//G_d(Q)$$ by coordinates $t_\omega$, where $\omega$ runs through representatives of the cyclic
equivalence classes of cycles. The quotient $M_d^{ssimp}(Q)$ parametrizes the closed orbits of $G_d$
in $R_d(Q)$, thus it parametrizes isomorphism classes of semisimple representations of $Q$ of dimension vector $d$.\\[1ex]
Formation of
the direct sum of representations induces direct sum maps
$$\bigoplus_{d,e}:M_d^{ssimp}(Q)\times M_e^{ssimp}(Q)\rightarrow M_{d+e}^{ssimp}(Q).$$
The complement of the image of all non-trivial direct sum maps is an open subset $$M_d^{simp}(Q)\subset M_d^{ssimp}(Q),$$
which, by definition, parametrizes isomorphism classes of simple representations of $Q$ of dimension vector $d$.
If $M_d^{simp}(Q)$ is non-empty, then, by the above, it is a smooth irreducible variety of dimension
$$\dim M_d^{simp}(Q)=\dim R_d(Q)-\dim PG_d=1-\langle d,d\rangle.$$
The variety $M_d^{simp}(Q)$ admits a model over ${\bf Z}$, that is, there exists a scheme over
${\bf Z}$ such that the base extension to ${\bf C}$ is isomorphic to $M_d^{simp}(Q)$. Its reduction
modulo prime powers therefore yields models of the moduli space $M_d^{simp}(Q)$ over finite fields,
whose rational points parametrize
isomorphism classes of absolutely simple representations of $Q$ of dimension vector $d$, where a representation is called
absolutely simple if it remains simple after scalar extension to an algebraic closure of the base field. The following
is proved in \cite{RCRP}:

\begin{theorem}\cite{RCRP} For all quivers $Q$ and all dimension vectors $d\in{\bf N}Q_0$, there exists a polynomial
$A_d(t)\in{\bf Z}[t]$ such that, for each finite field $k$, the number of $k$-rational points of $M_d^{simp}(Q)$ equals
$A_d(|k|)$. Moreover, the polynomials $A_d(t)$ are given by an explicit recursive formula.
\end{theorem}

To define a projectivized version of the moduli space $M_d^{simp}(Q)$, we have to recall some concepts of \cite{M}.
Given a representation of a reductive algebraic group $G$ on a complex vector space $V$, one has the notions
of stable, resp.~semistable, points for the induced action of $G$ on the projective space ${\bf P}V$, defined as follows:
a point $x\in{\bf P}V$ is called semistable if and only if for some (equivalently, for any) lift $v\in V$ of $x$,
the closure of the orbit $Gv$ does not contain $0\in V$. The point $x$ is called stable if it semistable and, additionally,
the orbit $Gv$ is closed and of maximal possible dimension. Denote by ${\bf P}V^{ss}$ (resp.~${\bf P}V^s$) the
set of semistable (resp.~stable) points of ${\bf P}V$. Defining the projective variety $Y$ by
$$Y:={\bf Proj}({\bf C}[V]^G),$$
there exists a quotient map $$\pi:{\bf P}V^{ss}\rightarrow Y,$$
and there exists a smooth open
subset $U$ of $Y$ such that the restriction of $\pi$ induces a geometric quotient $$\pi:{\bf P}V^s\rightarrow U.$$
We apply this construction of quotients to the projective space ${\bf P}R_d(Q)$ of the affine space $R_d(Q)$,
which inherits an action of $G_d$. The above mentioned results easily imply the following:

\begin{proposition} A point of ${\bf P}R_d(Q)$ is semistable if and only if some
(equivalently, any) lift to $R_d(Q)$ does not belong to
$N_d(Q)$. The stable points are precisely those
$x\in{\bf P}R_d(Q)$ whose lifts $X\in R_d(Q)$ correspond to simple representations.
\end{proposition}

\begin{definition} We define the projectivized moduli spaces
$${\bf P}M_d^{\rm simp}(Q)\mbox{ (resp.~}{\bf P}M_d^{\rm ssimp}(Q))$$
of simple (resp.~semisimple) representations of $Q$ of 
dimension vector $d$ as the
corresponding quotients of ${\bf P}R_d(Q)^{ss}$ (resp.~${\bf P}R_d(Q)^s$)
by the $G_d$-action, with quotient map $\pi_d(Q)$.
\end{definition}

The varieties ${\bf P}M_d^{ssimp}(Q)$ are introduced and analyzed in \cite{BS}
as moduli spaces of representations of $Q$ up to a notion of graded equivalence.\\[1ex]
Given a point $x\in{\bf P}M_d^{simp}(Q)$, we call a representation $X\in R_d(Q)$ a lift of $x$ if the corresponding point
$X\in{\bf P}R_d(Q)$ maps to $x$ via $\pi_d(Q)$.\\[1ex]
The quiver $Q$ is called strongly connected if and only if there exists a cycle through any given two vertices. A
particular example is the quiver consisting of a single cycle. The support ${\rm supp}(d)$ of a dimension vector
$d\in{\bf N}Q_0$ is defined as the full subquiver on vertices $i$ for which $d_i\not=0$.
As a direct consequence of the construction, a theorem of \cite{LBP1} implies the following criterion for non-emptyness
of ${\bf P}M_d^{simp}(Q)$:

\begin{theorem}\label{exsimples}\cite{LBP1}
The space ${\bf P}M_d^{\rm simp}(Q)$ is non-empty if and only if the following conditions are
satisfied:
\begin{itemize}
\item the support of $d$ is strongly connected, and contains at least one arrow,
\item either \begin{itemize}
\item $\langle i,d\rangle\leq 0$ and $\langle d,i\rangle\leq 0$ for all $i\in{\rm supp}(d)$, or
\item ${\rm supp}(d)$
is a single cycle and $d_i=1$ for all $i\in{\rm supp}(d)$.
\end{itemize}
\end{itemize}
\end{theorem}

The main aim of the paper is the following formula for the Euler characteristic in singular cohomology
with compact support of ${\bf P}M_d^{simp}(Q)$:

\begin{theorem}\label{euler} For all quivers $Q$ and all dimension vectors $d\in{\bf N}Q_0$, we
have
$$\chi_c({\bf P}M_d^{\rm simp}(Q))=|C_d^{prim}(Q)|,$$
the number of cyclic equivalence classes of primitive cycles
of dimension vector $d$.
\end{theorem}

As already mentioned in the introduction, this theorem will be proved by localization techniques. The key
ingredient is to describe the fixed point set of a natural torus action on ${\bf P}M_d^{simp}(Q)$, which will now be
described.\\[1ex]
We consider the torus $$T_Q:=({\bf C}^{*})^{|Q_1|},$$
whose elements will be written as tuples
$$t_\ast=(t_\alpha)_{\alpha\in Q_1}$$ parametrized
by the arrows in $Q$. Let $X(T_Q)$ be the group of characters of $T_Q$, which will always be identified with ${\bf N}Q_1$
via the identification of $\alpha\in{\bf N}Q_1$ with the character $$t_\ast\mapsto t_\alpha.$$
We have a natural right action of $T_Q$ on $R_d(Q)$ by rescaling the linear maps representing the arrows,
that is, $$(X_\alpha).t_\ast:=(t_\alpha\cdot X_\alpha)_\alpha.$$
Since this action obviously commutes with the $G_d$-action,
it descends to a right action of $T_Q$ on ${\bf P}M_d^{\rm simp}(Q)$.\\[1ex]
We call an element $\nu\in{\bf N}Q_1$ indivisible if $${\rm gcd}(\nu_\alpha\, :\, \alpha\in Q_1)=1.$$
Given such an indivisible element $\nu$, define a quiver $\widehat{Q}_\nu$ as follows: for the set of vertices, we have
$$(\widehat{Q}_\nu)_0=Q_0\times({\bf Z}Q_1/{\bf Z}\nu).$$
The arrows in $\widehat{Q}_\nu$ are given as
$$(\alpha,\overline{\lambda}):(i,\overline{\lambda})\rightarrow(j,\overline{\lambda+\alpha}),$$
where $(\alpha:i\rightarrow j)\in Q_1$ and $\overline{\lambda}\in{\bf Z}Q_1/{\bf Z}\nu$. We call this quiver the
almost universal abelian covering quiver of $Q$ associated to $\nu$ (omission of the factorization by ${\bf Z}\nu$ in the
definition constructs the universal abelian covering of $Q$). We have a natural action of
the abelian group ${\bf Z}Q_1$ on $(\widehat{Q}_\mu)_0$ given by
$$\mu\cdot(i,\overline{\lambda})=(i,\overline{\lambda+\mu}).$$
We consider dimension vectors
$$\widehat{d}=\sum_{i,\overline{\lambda}}\widehat{d}_{i,\overline{\lambda}}(i,\overline{\lambda})
\in{\bf N}(\widehat{Q}_\nu)_0$$
for $\widehat{Q}_\nu$ which lift $d$, by which we mean that 
$$\sum_{\overline{\lambda}}\widehat{d}_{i,\overline{\lambda}}=d_i$$
for each $i\in Q_0$. The action of ${\bf Z}Q_1$ on $(\widehat{Q}_\nu)_0$ induces an action on
${\bf N}(\widehat{Q}_\nu)_0$;
we call two dimension vectors equivalent if they belong to the same orbit under this action.

\begin{theorem}\label{fmr} The fixed point set
$${\bf P}M_d^{simp}(Q)^{T_Q}$$
is isomorphic to the disjoint union of moduli
$$\bigcup_{\nu,\widehat{d}}{\bf P}M_{\widehat{d}}^{simp}(\widehat{Q}_\nu),$$
where $\nu$ ranges over all indivisible elements of ${\bf N}Q_1$, and $\widehat{d}$ ranges over all equivalence classes
of dimension vectors for $\widehat{Q}_\nu$ lifting $d$.
\end{theorem}

\section{Detection of fixed points}\label{fps}

In this section, we analyze the closed subset
$${\bf P}M_d^{simp}(Q)^{T_Q}$$
of $T_Q$-fixed points in ${\bf P}M_d^{\rm simp}(Q)$ and show that they can be viewed
as representations of some covering quiver $\widehat{Q}_\nu$ for a dimension vector lifting $d$.\\[2ex]
The starting point for this analysis is the following lemma, which is inspired by \cite[Lemma 1.3]{Gi}:

\begin{lemma}\label{keylemma} Assume that $x$ is a $T_Q$-fixed point in ${\bf P}M_d^{\rm simp}(Q)$,
and let $X$ be a lift of $x$ in $R_d(Q)$.
Then there exists a tuple
$$((N_\alpha)_{\alpha\in Q_1},\psi,\varphi),$$
consisting of a tuple of non-negative integers $N_\ast=(N_\alpha)_{\alpha\in Q_1}$,
a map of algebraic groups $\varphi:T_Q\rightarrow G_d$, and a character $\psi\in X(T_Q)$,
such that for all $t_\ast\in T_Q$, we have
\begin{equation}X.(t_\alpha^{N_\alpha})=\psi(t_\ast)\cdot\varphi(t_\ast)X.\label{pe}\end{equation}
\end{lemma}

\proof Consider the subgroup
$$H\subset (G_d\times{\bf C}^*)\times T_Q$$ consisting of all pairs $((g,\lambda),t_\ast)$
such that
$$X.t_\ast=\lambda\cdot gX.$$
By the assumption of $x$ being a fixed point, the second projection
$p_2:H\rightarrow T_Q$ is surjective. Since $T_Q$, being a torus, is connected, commutative, and consists
of semisimple elements,
the successive restrictions of $p_2$ to the unit component of $H$, to its reductive part, and to its connected
center, all have to remain surjective.\\[1ex]
Therefore, there exists a torus $T$ in $H$ such that the restriction of
$p_2$ to $T$ remains surjective. 
But any surjection of tori is quasi-split, that is, there exist non-negative integers $N_\alpha$ for $\alpha\in Q_1$
and a map $s:T_Q\rightarrow T$ such that $$(p_2\circ s)(t_\ast)=(t_\alpha^{N_\alpha})$$
for all $t_\ast\in T_Q$. Composition of $s$ with the
projection $p_1:H\rightarrow G_d\times{\bf C}^*$ yields a map $$p_1s=(\varphi,\psi):T_Q\rightarrow G_d\times{\bf C}^*.$$
It is now an easy verification that $N_\ast=(N_\alpha)_\alpha$, $\psi$ and $\varphi$ fulfill the claimed property.\hb

We call a datum
$$(N_\ast=(N_\alpha)_\alpha,\psi\in X(T_Q),\varphi:T_Q\rightarrow G_d)$$ associated to a $T_Q$-fixed point
$x\in{\bf P}M_d^{simp}(Q)$ as in the above lemma (that is, fulfilling equation (\ref{pe})
for some lift $X\in R_d(Q)$ of $x$) a fixing datum associated to $x$.\\[2ex]
Now assume that a fixing datum $(N_\ast,\psi,\varphi)$ associated to $x$ is given, fulfilling equation (\ref{pe})
for a specified lift $X\in R_d(Q)$. We decompose $\varphi$ into
components
$$\varphi=(\varphi_i:T_Q\rightarrow{\rm GL}(V_i)).$$
Assume $X$ is given as a tuple of linear maps
$$(X_\alpha:V_i\rightarrow V_j)_{\alpha:i\rightarrow j}.$$
Then, by definition of the actions of $T_Q$ and $G_d$ on $R_d(Q)$, equation (\ref{pe}) can be reformulated as
\begin{equation}t_\alpha^{N_\alpha}\cdot X_\alpha=\psi(t_\ast)\cdot \varphi_j(t_\ast)X_\alpha\varphi_i(t_\ast)^{-1}
\label{eqn1}\end{equation}
for all $(\alpha:i\rightarrow j)\in Q_1$ and all $t_\ast\in T_Q$.\\[1ex]
The map $\varphi_i$ provides the vector space $V_i$ with the structure of a representation of $T_Q$,
for which we have the weight space decomposition
$$V_i=\bigoplus_{\lambda\in X(T_Q)}V_i^\lambda,$$
where $V_i^\lambda$ denotes the
$\lambda$-weight space of $V_i$ for each $\lambda\in X(T_Q)$. Formula (\ref{eqn1}) can now be used to detect a kind of
grading on the representation $X$, given by the above $\lambda$-weight spaces:

\begin{lemma}\label{lwsd} For all arrows $\alpha:i\rightarrow j$ and all $\lambda\in X(T_Q)$, we have
\begin{equation}X_\alpha(V_i^\lambda)\subset V_j^{\lambda+N_\alpha\alpha-\psi}.\label{wsd}\end{equation}
\end{lemma}

\proof If $v\in V_i^\lambda$, then $\varphi_i(t_\ast)v=\lambda(t_\ast)\cdot v$ by definition. Thus, for all
$t_\ast\in T_Q$, we have
\begin{eqnarray*}
\varphi_j(t_\ast)X_\alpha v&=&\varphi_j(t_\ast)X_\alpha\varphi_i(t_\ast)^{-1}\varphi_i(t_\ast)v\\
&=&(\psi(t_\ast)^{-1}t_\alpha^{N_\alpha}\cdot X_\alpha)(\lambda(t_\ast)\cdot v)\\
&=&(\psi(t_\ast)^{-1}t_\alpha^{N_\alpha}\lambda(t_\ast))\cdot X_\alpha v,
\end{eqnarray*}
which yields the claim.\hb

The following lemma shows that the character $\psi$ in a fixing datum is already determined by the tuple $N_\ast$.

\begin{corollary}\label{nu} Given a fixing datum as above, there exists a unique primitive element
$$\nu=\sum_\alpha\nu_\alpha\alpha\in X(T_Q)$$ such that $$\psi=\frac{1}{|\nu|}\sum_\alpha\nu_\alpha N_\alpha\alpha.$$
\end{corollary}

\proof Since $X$ belongs to ${\bf P}R_d^{ss}(Q)$ by definition, there exists a cycle $$\omega=(\alpha_1,\ldots,\alpha_s)$$
at a vertex $i\in Q_0$ in $Q$ such that $${\rm tr}(X_{\alpha_s}\cdot\ldots\cdot X_{\alpha_1})\not=0.$$
Denote by $\nu\in{\bf Z}Q_1$ the weight of $\omega$, thus $|\nu|=s$. By iterated application of Lemma \ref{lwsd},
we have $$(X_{\alpha_s}\cdot\ldots\cdot X_{\alpha_1})(V_i^{\lambda})\subset V_i^{\mu}$$ for the weight $\mu$ defined by
$$\mu=\lambda+\sum_{\alpha}\nu_\alpha N_\alpha\alpha-s\psi.$$
By non-vanishing of the trace of this operator, we have equality of the two weights $\lambda$, $\mu$,
and thus
$$\psi=\frac{1}{|\nu|}\sum_\alpha\nu_\alpha N_\alpha\alpha.$$
Dividing $\nu$ by ${\rm gcd}(\nu_\alpha\, :\, \alpha\in Q_1)$ if neccessary yields the existence of the
desired element $\nu$. Uniqueness of such an element $\nu$ follows from its indivisibility.\hb

A crucial point in the description of the fixed points is the question to what extend
a fixing datum $(N_\ast,\psi,\varphi)$ attached to a fixed point $x$
and a chosen lift $X$ is unique.\\[1ex]
We first study the effect of changing the lift $X$. Assume $$X'=\lambda\cdot g_0X$$ is another lift, for a non-zero
scalar $\lambda\in{\bf C}^*$ and a group element $g_0\in G_d$. It follows immediately that equation
(\ref{pe}), considered for the lift $X'$, is fulfilled for the fixing datum
\begin{equation}(N_\ast,\psi,g_0\varphi g_0^{-1}).\label{equiv1}\end{equation}
Next, choose a tuple $m_\ast=(m_\alpha)_\alpha$ of non-zero integers and consider the isogeny $i:T_Q\rightarrow T_Q$
defined by
$$t_\ast\mapsto t_\ast^{m_\ast}:=(t_\alpha^{m_\alpha})_\alpha.$$
A direct calculation then shows that the tuple
\begin{equation}((N_\alpha m_\alpha)_\alpha,\psi\circ i,\varphi\circ i)\label{equiv2}\end{equation}
is also a fixing datum associated to $x$ for the lift $X$.\\[1ex]
Finally, we can use the fact that the action of
$G_d$ on $R_d(Q)$ factors through $PG_d$. This allows us to choose an arbitrary character $\mu\in X(T_Q)$ and to replace
the given fixing datum by
\begin{equation}(N_\ast,\psi,\mu\cdot\varphi)\label{equiv3}\end{equation}
for the same lift $X$ of $x$.\\[1ex]
This analysis motivates the following:
\begin{definition} Two fixing data associated to a fixed point $x$ are called equivalent if they are equivalent
under the equivalence relation generated by the above three
operations, that is, by replacing $(N_\ast,\psi,\varphi)$ by (\ref{equiv1}), or (\ref{equiv2}), or (\ref{equiv3}) above.
\end{definition}

Using this
notion of equivalence, we have:

\begin{lemma} Any two fixing data for a $T_Q$-fixed point $x\in{\bf P}M_d^{simp}(Q)$ are equivalent.
\end{lemma}

\proof Let $$(N_\ast,\psi,\varphi),\; (N_\ast',\psi',\varphi')$$ be fixing data for $x$ fulfilling equation
(\ref{pe}) for lifts $X$, resp.~$X'$, of $x$. Since the map $$\pi_d(Q):{\bf P}R_d(Q)^{s}\rightarrow{\bf P}M_d^{simp}(Q)$$
is a geometric quotient, we find a non-zero scalar $\lambda\in{\bf C}^*$ and a group element $g_0\in G_d$
such that $$X'=\lambda\cdot g_0X.$$
Thus, by operation (\ref{equiv1}) on fixing data, we can assume $X'=X$
without loss of generality. We can use operation (\ref{equiv2}) on fixing data by twisting
$(N_\ast,\psi,\varphi)$ with the isogeny $(t_\ast\mapsto t_\ast^{N_\ast'})$, and similarly twisting
$(N_\ast',\psi',\varphi*)$
with the isogeny $(t_\ast\mapsto t_\ast^{N_\ast})$. This allows us, again without loss of generality, to assume
$N_\ast=N_\ast'$. So we arrive at the following equations:
$$X.t_\ast^{N_\ast}=\psi(t_\ast)\cdot\varphi(t_\ast)X,\;\;\; X.t_\ast^{N_\ast}=\psi'(t_\ast)\cdot\varphi'(t_\ast)X,$$
giving
$$\frac{\psi(t_\ast)}{\psi'(t_\ast)}\cdot X=\varphi(t_\ast)^{-1}\varphi'(t_\ast)X.$$
Applying Corollary \ref{nu} to both $\psi$ and $\psi'$, we arrive at $\psi=\psi'$, and thus at
$$X=\varphi(t_\ast)^{-1}\varphi'(t_\ast)X.$$
But since $X$ is a simple representation, its automorphism group reduces
to the scalars, so $$\varphi'=\mu\cdot\varphi'$$
for some character $\mu\in X(T_Q)$. Via operation (\ref{equiv3}) on fixing data,
we conclude the desired equivalence.\hb

Define an auxilliary quiver $\Gamma$ as follows: its set of vertices is
$$\Gamma_0=Q_0\times X(T_Q),$$
and for each arrow $(\alpha:i\rightarrow j)\in Q_1$ and each $\lambda\in X(T_Q)$, we define an arrow
$$(\alpha,\lambda):(i,\lambda)\rightarrow (j,\lambda+N_\alpha\alpha-\psi)$$ in $\Gamma_1$.
Moreover, we define a dimension vector $\widehat{e}\in{\bf N}\Gamma$ by
$$\widehat{e}_{i,\lambda}:=\dim_{\bf C}V_i^{\lambda}.$$
Lemma \ref{lwsd} now can be reinterpreted as saying that
$X$ can be viewed as a representation of $\Gamma$ of dimension vector $\widehat{e}$,
which is neccessarily a simple representation of $\Gamma$ since $X$ is so.
We will now show that $X$ is already supported on a subquiver
$\Gamma'$ of $\Gamma$ which depends only on $\nu$, and is thus independent of the fixing datum.\\[1ex]
Consider the map
$$h:X(T_Q)\rightarrow X(T_Q),\;\;\; h(\alpha):=N_\alpha\alpha-\psi,$$
and recall the indivisible element $\nu\in{\bf N}Q_1\simeq X(T_Q)$ attached to the fixed point in Corollary \ref{nu}. 

\begin{lemma}\label{lpc} The map $h$ induces an embedding $X(T_Q)/{\bf Z}\nu\rightarrow X(T_Q)$.
\end{lemma}

\proof We have
$$h(\nu)=\sum_\alpha\nu_\alpha(N_\alpha\alpha-\psi)=\sum_\alpha\nu_\alpha N_\alpha\alpha-(\sum_\alpha\nu_\alpha)\psi=0$$
by Corollary \ref{nu}. Now assume that $h(\lambda)=0$ for a character $\lambda\in X(T_Q)$. Again by Corollary \ref{nu},
we have
$$\sum_\alpha\lambda_\alpha N_\alpha\alpha=|\lambda|\cdot\psi=
\sum_\alpha\frac{|\lambda|}{|\nu|}\nu_\alpha N_\alpha\alpha.$$
Since all $N_\alpha$ are non-zero, we have $$\lambda=\frac{|\lambda|}{|\nu|}\nu.$$
Since $\nu$ is indivisible by assumption, we conclude that $\lambda$ is an integer multiple
of $\nu$.\hb

We choose a vertex $i_0\in Q_0$ and a character $\lambda_0\in X(T_Q)$
such that
$$V_{i_0}^{\lambda_0}\not=0.$$
By Lemma \ref{lpc}, we see immediately that we have an embedding of quivers
$$H:\widehat{Q}_\nu\rightarrow\Gamma$$ by defining $H$ on vertices and on arrows via
$$H_0(i,\overline{\lambda})=(i,\lambda_0+h(\alpha)),\;\;\; H_1(\alpha,\overline{\lambda})=(\alpha,h(\lambda)).$$
The image $\Gamma'$ of $H$ is obviously connected,
since $\widehat{Q}_\nu$ is so. Since $X$ is a simple representation, its support is connected and meets $\Gamma'$ by
definition. We conclude that the support of $X$ is contained in $\Gamma'$. Thus,
$X$ can be viewed as a representation of the covering quiver $\widehat{Q}_\nu$,
whose dimension vector $\widehat{d}$ is given by
$$\widehat{d}_{i,\overline{\lambda}}=\widehat{e}_{i,\lambda_0+h(\lambda)}.$$
Note again that the choice of $\nu$ depends only on the given fixed point $x$ by Corollary \ref{nu}.\\[1ex]
We finally prove
that the dimension vector $\widehat{d}$ is determined uniquely up to the notion of equivalence of dimension vectors for
$\widehat{Q}_\nu$ introduced at the end of section \ref{projectivized}, for which we have to study the effect of the
operations (\ref{equiv1}), (\ref{equiv2}), (\ref{equiv3}) on fixing data on the dimension vector $\widehat{d}$:\\[1ex]
Obviously,
operation (\ref{equiv1}) does neither affect the weights $\lambda$ for which $V_i^{\lambda}\not=0$, nor the dimensions of
the $V_i^{\lambda}$. The twist by an isogeny in operation (\ref{equiv2}) is compensated for by the map $h$.
So it remains to
consider operation (\ref{equiv3}): it affects the weights occuring in the weight space decomposition of all $V_i$ by
translation by a fixed character $\mu$, and thus affects the dimension vector $\widehat{e}$ by translation by $\mu$.
Thus, the same holds for the dimension vector $\widehat{d}$, which precisely defines the notion of equivalence
on such dimension vectors.\\[2ex]
We summarize our present analysis of $T_Q$-fixed points
in ${\bf P}M_d^{simp}(Q)$ by:

\begin{proposition}\label{detect}
For each fixed point $x\in{\bf P}M_d^{simp}(Q)^{T_Q}$ and each lift $X$ of $x$ to $R_d(Q)$, there
exists a unique indivisible element $\nu\in {\bf N}Q_1$ and a unique (up to equivalence) dimension vector $\widehat{d}$ for
$\widehat{Q}_\nu$ lifting $d$,
such that $X$ can be viewed as a simple representation of $\widehat{Q}_\nu$ of dimension vector $\widehat{d}$.
\end{proposition} 

In the following section, we will formalize this way of viewing $X\in R_d(Q)$ as a representation of a covering quiver,
by embedding its projectivized moduli of simple representations into the $T_Q$-fixed points.

\section{Construction of fixed points}\label{construction}

Fix an indivisible element $\nu\in{\bf N}Q_1\simeq X(T_Q)$ as above,
and let $\widehat{d}$ be a dimension vector for $\widehat{Q}_\nu$
lifting a dimension vector $d$. Our aim is to construct a closed embedding
$$P_\nu:{\bf P}M_{\widehat{d}}^{simp}(\widehat{Q}_\nu)\rightarrow{\bf P}M_d^{simp}(Q)^{T_Q}.$$

We first construct a map on the level of varieties of representations.
Choose vector spaces $V_{i,\overline{\lambda}}$ of dimension $\widehat{d}_{i,\overline{\lambda}}$ for all $i\in Q_0$
and all $\overline{\lambda}\in X(T_Q)/{\bf Z}\nu$, and define
$$V_i:=\bigoplus_{\overline{\lambda}}V_{i,\overline{\lambda}}.$$
Given a representation
$$(X_{\alpha,\overline{\lambda}}:
V_{i,\overline{\lambda}}\rightarrow V_{j,\overline{\lambda+\alpha}})_{\alpha,\overline{\lambda}}$$ of
$\widehat{Q}_\nu$ of dimension vector $\widehat{d}$, we can define a representation of $Q$ of dimension vector $d$
by the linear maps $$(\bigoplus_{\overline{\lambda}}X_{\alpha,\overline{\lambda}}:V_i\rightarrow V_j)_{(\alpha:
i\rightarrow j)\in Q_1}.$$
This induces a linear, and thus regular and closed, map of representation varieties
$$P_\nu:R_{\widehat{d}}(\widehat{Q}_\nu)\rightarrow R_d(Q).$$
Using the given decomposition of the vector spaces $V_i$, we also have a natural embedding of $G_{\widehat{d}}$ into
$G_d$ via
$$(G_{\widehat{d}})_i=\prod_{\overline{\lambda}}{\rm GL}(V_{i,\overline{\lambda}})\rightarrow{\rm GL}(V_i)=(G_d)_i,$$
such that the map $P_\nu$ becomes equivariant for the induced actions of $G_{\widehat{d}}$.
Moreover, $P_\nu$ is obviously equivariant for the ${\bf C}^*$-action
simultaneously rescaling the arrows in $\widehat{Q}_\nu$ and in $Q$, respectively. By the universal property of
quotients, we thus get an induced closed map
$${P}_\nu:{\bf P}M_{\widehat{d}}^{ssimp}(\widehat{Q}_\nu)\rightarrow{\bf P}M_d^{ssimp}(Q).$$

We easily get a converse to Lemma \ref{keylemma}:

\begin{lemma} The image of $P_\nu$ is contained in ${\bf P}M_d^{ssimp}(Q)^{T_Q}$.
\end{lemma}

\proof Choose a tuple $N_\ast=(N_\alpha)_{\alpha\in Q_1}$ of non-negative integers such that
$$\psi:=\frac{1}{|\nu|}\sum_\alpha
\nu_\alpha N_\alpha\alpha\in{\bf Q}Q_1$$
already belongs to ${\bf N}Q_1$
(for example, choose $N_\alpha=|\nu|$ for all $\alpha\in Q_1$).
Define a map
$$\varphi=(\varphi_i:T_Q\rightarrow{\rm GL}(V_i))_{i\in Q_0}:T_Q\rightarrow G_d$$
by sending $t_\ast$ via $\varphi_i$ to the automorphism of $V_i$
which acts as multiplication by $h(\overline{\lambda})(t_\ast)$ on each $V_{i,\overline{\lambda}}$. Equivalently,
$\varphi_i$ is defined in such a way that each $V_{i,\overline{\lambda}}$ is precisely the $h(\overline{\lambda})$-weight
space in the $T_Q$-representation on $V_i$ defined by $\varphi_i$, converse to the construction preceding
Lemma \ref{lwsd}.\\[1ex]
By the properties of $h$ established in Lemma \ref{lpc}, this makes $\varphi$ well-defined.
Then we can easily verify that equation (\ref{pe}) holds for each representation $Y=P_\nu(X)$ in the image of
$P_\nu$. We have thus constructed a fixing datum for the image $y\in{\bf P}M_d^{simp}(Q)$ of $Y$. By the definition of
${\bf P}M_d(T_Q)$ as a categorical quotient, we thus have
$$y.t_\ast^{N_\ast}=y$$
for all $t_\ast\in T_Q$. Thus, $y$ belongs to the set of $T_Q$-fixed points.\hb

The following two propositions show that the above map $P_\nu$ induces an embedding on the projectivized moduli of
simple representations.

\begin{proposition}\label{simplicity}
The map $P_\nu$ maps ${\bf P}M_{\widehat{d}}^{simp}(\widehat{Q}_\nu)$ to ${\bf P}M_d^{simp}(Q)$.
\end{proposition}

\proof Let $X$ be a simple representation of
$\widehat{Q}_\nu$ and assume that $$U=(U_i)_{i\in Q_0}$$ is a non-zero simple subrepresentation of
dimension vector $e\in{\bf N}Q_0$ of $Y=P_\nu(X)$.
Since $Y$ lies in the set of $T_Q$-fixed points by the previous lemma, we can find a 
fixing datum $(N_\ast,\psi,\varphi)$ such that $Y$ fulfills equation (\ref{eqn1}), that is, we have
$$t_\alpha^{N_\alpha}\cdot Y_\alpha=\psi(t_\ast)\cdot \varphi_j(t_\ast)Y_\alpha\varphi_i(t_\ast)^{-1}$$
for each arrow $(\alpha:i\rightarrow j)\in Q_1$ and each $t_\ast\in T_Q$.
Fixing an element $t_\ast\in T_Q$, consider the
tuple of subspaces
$$\varphi(t_\ast)U=(\varphi_i(t_\ast)U_i)_{i\in Q_0}.$$
Using equation (\ref{pe}), we have
\begin{eqnarray*}
Y_\alpha\varphi_i(t_\ast)U_i&=&t_\alpha^{-N_\alpha}\psi(t_\ast)\cdot
\varphi_j(t_\ast)Y_\alpha\varphi_i(t_\ast)^{-1}\varphi_i(t_\ast)
U_i\\
&\subset&\varphi_j(t_\ast)Y_\alpha U_i\\
&\subset&\varphi_j(t_\ast)U_j,
\end{eqnarray*}
thus $\varphi(t_\ast)U$ is again a simple subrepresentation of $Y$.\\[1ex]
The same computation also yields the following: if $u$ denotes the point of ${\bf P}M_e^{simp}(Q)$
corresponding to $U$, the point $u.t_\ast^{N_\ast}$ corresponds to the representation $\varphi(t_\ast)U$. We can thus
conclude that $Y$ contains a direct sum of simple representations corresponding to all points
in the $T_Q$-orbit $u.T_Q$ of $u$ in ${\bf P}M_e^{simp}(Q)$ as a subrepresentation. But $Y$ is a finite dimensional
representation, and thus the orbit $u.T_Q$ is finite. It is also connected, since $T_Q$ is so,
and we can conclude that $u$
is a $T_Q$-fixed point in ${\bf P}M_e^{simp}(Q)$. This in turns yields $$\varphi(t_\ast)U=U$$ for all
$t_\ast\in T_Q$.\\[1ex]
Thus, we can regard each $U_i\subset V_i$ as a subrepresentation of the representation of $T_Q$ on $V_i$ defined by
$\varphi_i$, and hence $U_i$ is compatible with the corresponding decomposition into isotypical components.
Since the map $h$ is injective by Lemma \ref{lpc}, the spaces $V_{i,\overline{\lambda}}$ are precisely these
isotypical components, and thus $U_i$ is compatible with the $V_{i,\overline{\lambda}}$, that is,
$$U_i=\bigoplus_{\overline{\lambda}}(U_i\cap V_{i,\overline{\lambda}}).$$
We conclude that $U$ can already by regarded as a non-zero subrepresentation of $X$, and thus $U=X$ by simplicity
of $X$. This proves simplicity of $Y$.\hb

\begin{proposition}\label{injectivity}
If $X$ and $X'$ are simple representations of $\widehat{Q}_\nu$ such that the $Q$-representations
$P_\nu(X)$ and $P_\nu(X')$ are
isomorphic, then $X$ and $X'$ are already isomorphic.
\end{proposition}

\proof The strategy of this proof is similar to the previous one. Denote $Y=P_\nu(X)$ and $Y'=P_\nu(X')$. Let
$$\gamma=(\gamma_i\in{\rm GL}(V_i))_{i\in Q_0}$$
be an isomorphism of $Y$ and $Y'$. By definition, this means
$$Y'_\alpha\circ\gamma_i=\gamma_j\circ Y_\alpha$$
for all $(\alpha:i\rightarrow j)\in Q_1$. Choose a fixing datum $(N_\ast,\psi,\varphi)$ such that $Y$ and $Y'$ fulfill
equation (\ref{eqn1}). Then, for all $t_\ast\in T_Q$ and all $(\alpha:i\rightarrow j)\in Q_1$, we have:
\begin{eqnarray*}
\psi(t_\ast)\cdot \varphi_j(t_\ast)Y'_\alpha\varphi_i(t_\ast)^{-1}\gamma_i&=&t_\alpha^{N_\alpha}\cdot Y'_\alpha\gamma_i\\
&=&\gamma_j
(t_\alpha^{N_\alpha}\cdot Y_\alpha)\\
&=&\psi(t_\ast)\cdot\gamma_j\varphi_j(t_\ast)Y_\alpha\varphi_i(t_\ast)^{-1},
\end{eqnarray*}
which gives
$$Y'_\alpha(\varphi_i(t_\ast)^{-1}\gamma_i\varphi_i(t_\ast))=(\varphi_j(t_\ast)^{-1}\gamma_j\varphi_j(t_\ast))Y_\alpha$$
for all $\alpha\in Q_1$. Since both $Y$ and $Y'$ are simple representations by Proposition \ref{simplicity},
we can apply
Schur's Lemma to conclude that
$$\varphi(t_\ast)^{-1}\gamma\varphi(t_\ast)=(\varphi_i(t_\ast)^{-1}\gamma_i\varphi_i(t_\ast))_{i\in Q_0}$$
is a scalar multiple of $\gamma$ for all $t_\ast\in T_Q$. Similar to the proof of Lemma \ref{keylemma},
we consider the subgroup
$$T\subset T_Q\times{\bf C}^*$$
consisting of pairs $(t_\ast,a)$ such that
$$\varphi(t_\ast)^{-1}\gamma\varphi(t_\ast)=a\cdot \gamma.$$
By the above, the projection $p$ of $T$ to $T_Q$ is surjective. It is obviously injective, thus it defines an isomorphism
$$p:T\stackrel{\sim}{\rightarrow} T_Q.$$
Composing the inverse of $p$ with the projection of $T$ to ${\bf C}^*$ yields a character $a$
of $T_Q$ such that
$$\varphi(t_\ast)^{-1}\gamma\varphi(t_\ast)=a(t_\ast)\cdot\gamma$$
for all $t_\ast\in T_Q$. As in the proof of Proposition \ref{simplicity}, we now consider each vector space
$V_i$ as a representation of $T_Q$ via $\varphi_i$. It is now easy to see that each $\gamma_i$ induces an isomorphism
between the $\lambda$-weight space and the $(\lambda+a)$-weight space, for each $\lambda\in X(T_Q)$. But $V_i$ is
finite dimensional, thus $$a=0\in X(T_Q).$$
We conclude that the isomorphism $\gamma$ fixes each weight space, and can
thus be regarded as an isomorphism of the representations $X$ and $X'$ of $\widehat{Q}_\nu$.\hb

This proposition allows us to conclude that the map 
$$P_\nu:{\bf P}M_{\widehat{d}}^{simp}(\widehat{Q}_\nu)\rightarrow{\bf P}M_d^{simp}(Q)$$
is a closed embedding.\\[1ex]
We can now combine the results of this and the previous section to prove Theorem \ref{fmr}:

\begin{theorem} The fixed point set
$${\bf P}M_d^{simp}(Q)^{T_Q}$$
is isomorphic to the disjoint union of moduli
$$\bigcup_{\nu,\widehat{d}}{\bf P}M_{\widehat{d}}^{simp}(\widehat{Q}_\nu),$$
where $\nu$ ranges over all indivisible elements of ${\bf N}Q_1$, and $\widehat{d}$ ranges over all equivalence classes
of dimension vectors for $\widehat{Q}_\nu$ lifting $d$.
\end{theorem}

\proof By the previous proposition, we have a system of closed embeddings
$$(P_\nu:{\bf P}M_{\widehat{d}}^{simp}(\widehat{Q}_\nu)\rightarrow{\bf P}M_d^{simp}(Q))_{\nu,\widehat{d}},$$
where $\nu$ runs over all indivisible elements in ${\bf N}Q_1$, and $\widehat{d}$ runs
over all equivalence classes of dimension vectors $\widehat{d}$ for $\widehat{Q}_\nu$ lifting $d$.
By Proposition \ref{detect}, each fixed point $x\in{\bf P}M_d^{simp}(Q)^{T_Q}$ belongs to the image of a unique such
embedding. Thus, the images of the various $P_\nu$ cover
${\bf P}M_d^{simp}(Q)^{T_Q}$ and are pairwise disjoint. Thus, they are
exactly the connected components of ${\bf P}M_d^{simp}(Q)^{T_Q}$.\hb

\example Consider four-dimensional representations of
the quiver with one vertex and two loops $\alpha,\beta$. The criterion Theorem \ref{exsimples}
for non-emptyness of a moduli space
${\bf P}M_d^{simp}(Q)$ puts strong restrictions on the pairs $(\nu,\widehat{d})$ which have to be considered.
We get the following five non-empty fixed point components (an arrow `$\Rightarrow$' (resp.~`$\rightarrow$')
represents an arrow of
$\widehat{Q}_\nu$ covering $\alpha$ (resp.~$\beta$); the numbers at the vertices represent entries of the dimension
vectors):
\begin{description}
\item[$Q^1$:] $\begin{array}{ccccc}1&\stackrel{\Rightarrow}{\leftarrow}&2&\stackrel{\Rightarrow}{\leftarrow}&1
\end{array}$
\item[$Q^2$:] $\begin{array}{ccccccc}1&\stackrel{\Rightarrow}{\leftarrow}&1&\stackrel{\Rightarrow}{\leftarrow}&
1&\stackrel{\Rightarrow}{\leftarrow}&1\end{array}$
\item[$Q^3$:] $\begin{array}{ccccccc}&&1&&\Rightarrow&&1\\ &\swarrow&&\nwarrow&&\swarrow&\\
1&&\Rightarrow&&1&\end{array}$
\item[$Q^4$:] $\begin{array}{ccc}1&\Rightarrow&1\\ \uparrow&&\Downarrow\\ 1&\Leftarrow&1\end{array}$
\item[$Q^5$:] $\begin{array}{ccc}1&\Rightarrow&1\\ \uparrow&&\downarrow\\ 1&\leftarrow&1\end{array}$
\end{description}

\section{Computation of the Euler characteristic}\label{ind}

The aim of this section is to prove Theorem \ref{euler}, utilizing the localization theorem \ref{fmr}. Therefore, we have
to relate classes of primitive cycles in $Q$ and in the covering quivers $\widehat{Q}_\nu$.\\[1ex]
We start by studying cycles in the covering quivers $\widehat{Q}_\nu$.
We have a natural map of quivers
$$\Pi_\nu:\widehat{Q}_\nu\rightarrow Q$$
forgetting the class of the character $\lambda$; more precisely, $\Pi_\nu$ maps a vertex
$(i,\overline{\lambda})$ to the vertex $i$ and an arrow
$(\alpha,\overline{\lambda})$ to the arrow $\alpha$. Note that, by definition of $\widehat{Q}_\nu$,
fixing an arrow $\alpha:i\rightarrow j$ in $Q_1$ and a
vertex $(i,\overline{\lambda})\in(\widehat{Q}_\nu)_0$ lifting $i$, there exists a unique arrow
$$(\alpha,\overline{\lambda}):(i,\overline{\lambda})\rightarrow(j,\overline{\lambda+\alpha})$$
in $\widehat{Q}_\nu$ lifting $\alpha$ and starting at $(i,\overline{\lambda})$.\\[1ex]
Now assume that
$$\widehat{\omega}\, :\,(i_0,\overline{\lambda_0})
\stackrel{(\alpha_1,\overline{\lambda_0})}{\rightarrow}
(i_1,\overline{\lambda_1})
\stackrel{(\alpha_2,\overline{\lambda_1})}{\rightarrow}
\ldots
\stackrel{(\alpha_s,\overline{\lambda_{s-1}})}{\rightarrow}
(i_s,\overline{\lambda_s})=(i_0,\overline{\lambda_0})$$
is a cycle in $\widehat{Q}_\nu$. Then, by the definition of the arrows in $\widehat{Q}_\nu$, we have
$$\lambda_k=\lambda_0+\sum_{l=1}^k\alpha_l\in{\bf Z}Q_1/{\bf Z}\nu,$$
and thus
$$\sum_{l=1}^s\alpha_l\in{\bf N}\nu.$$
Denoting by $\omega=(\alpha_1,\ldots,\alpha_s)$ the projection of $\widehat{\omega}$ via $\Pi_\nu$,
we see that
$$|\omega|\in{\bf N}\nu.$$
To summarize, we have proved:

\begin{lemma}\label{cqn}
A cycle $\widehat{\omega}$ in $\widehat{Q}_\nu$ projecting to a cycle $\omega$ in $Q$ is uniquely determined
by $\omega$ and its initial vertex $(i_0,\overline{\lambda_0})\in{\bf N}(\widehat{Q}_\nu)_0$. Moreover, the weight
$|\omega|$ of $\omega$ is a multiple of $\omega$.
\end{lemma}

As in the previous section, we now consider the set of all pairs $(\nu,\widehat{d})$
consisting of an indivisible element in ${\bf N}Q_1$ and an equivalence class of a dimension vector $\widehat{d}$ for
$\widehat{Q}_\nu$ lifting $d$.

\begin{proposition}\label{procycles} There exists a bijection between $C_d^{prim}(Q)$ and
$$\bigcup_{\nu,\widehat{d}}C_{\widehat{d}}^{prim}(\widehat{Q}_\nu).$$
\end{proposition}

\proof Obviously, the projection $\Pi_\nu$ induces a well-defined map
$$p_{\nu,\widehat{d}}:C_{\widehat{d}}(\widehat{Q}_\nu)\rightarrow C_d(Q)$$
from cycle classes in $C_{\widehat{d}}(\widehat{Q}_\nu)$
to cycle classes in $C_d(Q)$.\\[1ex]
Now assume $\widehat{\omega}$
is a primitive cycle of length $s$ in $\widehat{Q}_\nu$.
By the previous lemma, the weight $|\omega|$ of $\omega=p_{\nu,\widehat{d}}(\widehat{\omega})$
is a multiple of $\nu$.
We claim that $\omega$ is again primitive. Assume to the contrary that
$$\omega=(\omega')^n$$
for a cycle $\omega'$ in $Q$ of length $s'$ and an integer $n>1$ such that $s=ns'$.
In particular, the weight $|\omega'|$ of $\omega'$ divides $|\omega|$,
and thus $|\omega'|$ is also a multiple of $\nu$, since $\nu$ in indivisible. Using Lemma \ref{cqn} again, it follows
that there exists a unique lift $\widehat{\omega}'$ of $\omega'$ to $\widehat{Q}_\nu$
having the same initial vertex as $\widehat{\omega}$. The $n$-th power of $\widehat{\omega}'$ therefore equals
$\widehat{\omega}$, contradicting primitivity.\\[1ex]
So each $p_{\nu,\widehat{d}}$ restricts to a map
$$p_{\nu,\widehat{d}}:C_{\widehat{d}}^{prim}(\widehat{Q}_\nu)\rightarrow C_d^{prim}(Q).$$
To prove the claimed bijection, we will now show that each class of a primitive cycle
$\omega=(\alpha_1,\ldots,\alpha_s)$ at $i\in Q_0$ in $Q$ admits a lift via
a unique (up to equivalence of dimension vectors $\widehat{d}$) map $p_{\nu,\widehat{d}}$.\\[1ex]
To prove existence, we consider
the unique primitive element $\nu$ such that $|\omega|$ is a multiple of $\nu$, and we choose
an arbitrary $\lambda\in{\bf Z}Q_1$. Lemma \ref{cqn} guarentees existence of a unique cycle in $\widehat{Q}_\nu$
lifting $\omega$ and having inital vertex $(i,\overline{\lambda})$, which is neccessarily primitive, since $\omega$
is so. Its dimension
vector is obviously a lift of $d$.\\[1ex]
Conversely, assume that the class of $\omega$ belongs to the
image of some $p_{\nu,\widehat{d}}$. Then $|\omega|$ is a multiple of $\nu$, defining $\nu$ uniquely.
A final application of
Lemma \ref{cqn} shows that a lifting cycle is uniquely determined
up to some translation in the dimension vector $\widehat{d}$.\hb

We need another preliminary lemma, which will function as the starting point of an induction.

\begin{lemma}\label{indanf} Assume that
$${\bf P}M_d^{simp}(Q)\not=\emptyset,$$
and that there exists a lift $\widehat{d}$ of $d$ to some $\widehat{Q}_\nu$ such that
$$|{\rm supp}(\widehat{d})|=|{\rm supp}(d)|.$$
Then ${\rm supp}(d)$ is a single cycle quiver, and $d_i=1$ for all $i$ in the support of $d$.
\end{lemma}

\proof Assume without loss of generality that the support of $d$ equals $Q$.
Equality of the cardinality of the supports implies that, for each $i\in Q_0$, there exists a unique
$\overline{\lambda}(i)\in{\bf Z}Q_1/{\bf Z}\nu$ such that $\widehat{d}_{i,\overline{\lambda}(i)}\not=0$. By definition
of the arrows in $\widehat{Q}_\nu$, we see that the projection $\Pi_\nu$ induces an isomorphism of quivers
$$\Pi_\nu:{\rm supp}(\widehat{d})\stackrel{\sim}{\rightarrow}Q.$$
Since there exists a simple representation of
$Q$ of dimension vector $d$, the quiver $Q$ is strongly connected by Theorem \ref{exsimples}. Choosing a cycle $\omega$
in $Q$, Lemma \ref{cqn} and the above isomorphism of quivers show that $|\omega|$ is already a multiple of
$\nu\in{\bf N}Q_1$. Thus, any two cycles in $Q$ have proportional weights. But this already implies that
$Q$ is a single cycle quiver. Theorem \ref{exsimples} now provides the claimed properties of $d$.\hb

After these combinatorial preparations, we turn to the topological aspects.\\[1ex]
Proofs of the following localization principle can be found, for example, in \cite{CG}, \cite{EM}:

\begin{theorem}\label{loca} Let $X$ be a complex variety on which a torus $T$ acts.
Then, for the Euler characteristic in cohomology
with compact support, we have the equality
$$\chi_c(X)=\chi_c(X^T).$$
\end{theorem}

Moreover, by the long exact sequence for cohomology with compact support,
the Euler characteristic behaves additively with repect to disjoint unions. Theorem \ref{fmr} therefore immediately
implies:

\begin{equation}\chi_c({\bf P}M_d^{simp}(Q))=\sum_{\nu,\widehat{d}}\chi_c({\bf P}M_{\widehat{d}}^{simp}(\widehat{Q}_\nu)),
\label{eqeuler}\end{equation}
where $\nu$ ranges over all primitive elements in ${\bf N}Q_1$, and $\widehat{d}$ ranges over all equivalence classes
of dimension vectors for $\widehat{Q}_\nu$ lifting $d$.\\[1ex]
This formula is the final ingredient for the proof of Theorem \ref{euler}:

\begin{theorem} For all quivers $Q$ and all dimension vectors $d\in{\bf N}Q_0$, we have
$$\chi_c({\bf P}M_d^{simp}(Q))=|C_d^{prim}(Q)|.$$
\end{theorem}

\proof We fix a nonnegative integer $b$ and prove the theorem for all quivers $Q$ and all dimension vectors $d$ for $Q$
such that $|d|=b$. We proceed by downward induction on $|{\rm supp}(d)|$.
The induction therefore starts with the case
$$|{\rm supp}(d)|=b,$$
which means that
$d_i=1$
for all vertices
$i$ in the support of $d$. A primitive cycle of dimension vector $d$ for $Q$ is thus a cycle passing through
each vertex of $Q$ precisely once. Such a cycle exists if and only if $Q$ is a single cycle quiver. Thus we have
$$|C_d^{prim}(Q)|=\left\{\begin{array}{ccl}1&,&\mbox{$Q$ is a single cycle quiver,}\\
0&,&\mbox{otherwise.}\end{array}\right\}$$
Assume that $\chi_c({\bf P}M_d^{simp}(Q))$ is non-zero. By Theorem \ref{loca}, its set of $T_Q$-fixed points is non-zero,
too, and Theorem \ref{fmr} implies that there exists a lift $\widehat{d}$ of $d$ to some covering quiver
$\widehat{Q}_\nu$. Since $d_i=1$ for all $i$ in the support of $d$, we have
$$|{\rm supp}(\widehat{d})|=|{\rm supp}(d)|,$$
and we can apply Lemma \ref{indanf} to conclude that $Q$ is a single cycle quiver. In this case,
the moduli space
reduces to a point (for example, by the dimension formula for the moduli space in section \ref{projectivized}),
thus the Euler characteristic equals $1$. Therefore, the claimed equality holds.\\[1ex]
For the inductive step, assume now that $|{\rm supp}(d)|$ is strictly less than $b$, and compute the Euler
characteristic by Formula \ref{eqeuler}. Assume that a pair $(\nu,\widehat{d})$ contributes to the sum in this formula.
Applying Lemma \ref{indanf} again, we see that
$$|{\rm supp}(\widehat{d})|>|{\rm supp}(d)|,$$
and thus the inductive assumption applies to give
$$\chi_c({\bf P}M_d^{simp}(Q))=\sum_{\nu,\widehat{d}}|C_{\widehat{d}}^{prim}(\widehat{Q}_\nu)|.$$
Using Proposition \ref{procycles}, the claim follows.\hb

\example In continuation of the example of the previous section,
we see that the second and the third components have Euler characteristic
equal to $0$. The fourth and fifth components are single cycles, so they both contribute by $1$ to the overall Euler
characteristic. The first component has to be further localized, yielding a four-cycle. Thus, we see that
the Euler characteristic $\chi_c({\bf P}M_d^{simp}(Q))$ equals $3$.

\section{Comments}

In this final section, we comment on the two main results Theorem \ref{fmr} and Theorem \ref{euler} and give some
examples.\\[2ex]
The first comment concerns a more structural interpretation of Theorem \ref{euler}: the $0$-th Hochschild homology
of the path algebra of a quiver $Q$ can be computed as the path algebra modulo the linear subspace spanned by commutators.
It is easy to see that this space is spanned by cyclic equivalence classes of cycles in $Q$, thus
$${\rm HH}_0({\bf C}Q)={\bf C}Q/[{\bf C}Q,{\bf C}Q]=\bigoplus_{d\in{\bf N}Q_0}\bigoplus_{\overline{\omega}\in C_d(Q)}
{\bf C}\omega,$$
yielding a natural ${\bf N}Q_0$-grading on ${\rm HH}_0({\bf C}Q)$.\\[1ex]
The prominent role of the primitive cycles in the present setup then draws attention to a (rather weak)
additional structure
on ${\rm HH}_0({\bf C}Q)$, namely the $p$-th power maps $T_p$, mapping the cyclic equivalence class
of a cycle $\omega$ to the class of $\omega^p$. This operation allows us to define the primitive part
$${\rm HH}_0({\bf C}Q)^{prim}$$
as the subspace spanned by primitive cycles, and Theorem \ref{euler} can be rewritten in more algebraic terms as
$$\chi_c({\bf P}M_d^{prim}(Q))=\dim {\rm HH}_0({\bf C}Q)^{prim}_d.$$

\example In the case of the $m$-loop quiver $Q_m$, the Euler characteristic $\chi_c({\bf P}M_d^{simp}(Q_m))$
equals the number of primitive necklaces of length $d$ with beads of $m$ colours, a number which appears in many
algebraic and combinatorial interpretations. One prominent such interpretation is
the dimension of the degree $d$-component of the free Lie algebra in $m$ generators
(see \cite{Reut}). We have the explicit formula:
$$\chi_c({\bf P}M_d^{simp}(Q_m))=\frac{1}{d}\sum_{r|d}\mu(\frac{d}{r})m^r,$$
where $\mu$ denotes the number-theoretic Moebius function.\\[1ex]
Given this relation with the free Lie algebra $L^{(m)}$,
one can speculate whether $L^{(m)}$ or some related algebra might be constructed using some kind of convolution
construction on the moduli spaces ${\bf P}M_d^{simp}(Q_m)$.\\[2ex]
Non-emptyness of ${\bf P}M_d^{simp}(Q)$ does not neccessarily imply non-vanishing of the Euler characteristic,
as the following example shows:\\[1ex]
\example Consider the quiver $Q$ defined as
$$\begin{array}{ccccc}&\alpha&&\beta&\\ i&\stackrel{\rightarrow}{\leftarrow}&j&\stackrel{\rightarrow}{\leftarrow}&
k\\ &\delta&&\gamma& \end{array}$$
and the dimension vector $d=i+j+k$. The quotient map
$${\bf P}R_d(Q)\rightarrow {\bf P}M_d^{ssimp}(Q)$$
is given by
$${\bf P}^3\rightarrow{\bf P}^1,\;\;\; (X_\alpha:X_\beta:X_\gamma:X_\delta)\mapsto(X_\alpha X_\delta:X_\beta
X_\gamma).$$
Simplicity of a representation
is equivalent to non-vanishing of the four scalars representing the arrows, thus
$${\bf P}M_d^{simp}(Q)\simeq{\bf P}^1\setminus\{0,\infty\},$$
and we see that the Euler characteristic equals $0$.\\[2ex]
Given the description of the torus fixed point set, it is natural to ask whether it can be used to find more detailed
information on the cohomology of ${\bf P}M_d^{simp}(Q)$ than just the Euler
characteristic. A natural approach is to
use localization theorems in equivariant cohomology (see for example \cite{GKM}).
For these to be applicable, and in particular for relating the results
to non-equivariant cohomology, a key property is equivariant formality of ${\bf P}M_d^{simp}(Q)$, as defined in
\cite{GKM}. This property means that the $T_Q$-equivariant cohomology
$$H_{T_Q}^\ast({\bf P}M_d^{simp}(Q))$$ is a free module over the $T_Q$-equivariant cohomology ring of a point
$H_{T_Q}^\ast(pt)$. It implies the degeneration of the natural spectral sequence for equivariant cohomology, yielding
$$H_{T_Q}^\ast({\bf P}M_d^{simp}(Q))\simeq H^\ast({\bf P}M_d^{simp}(Q))\otimes H_{T_Q}^\ast(pt).$$
In continuation of the above example, we see that equivariant formality does not hold in the present situation:\\[1ex]
\example The torus $T_Q$ acts transitively on the quotient ${\bf P}^1\setminus\{0,\infty\}$ by
$$(x:y).(t_\alpha,t_\beta,t_\gamma,t_\delta)=(t_\alpha t_\delta x:t_\beta t_\gamma y),$$
with stabilizer of $1\in{\bf P}^1$ being the subtorus
$$T:=\{(t_\alpha,t_\beta,t_\gamma,t_\delta)\, :\, t_\alpha t_\delta=t_\beta t_\gamma\}.$$
So the $T_Q$-equivariant cohomology of ${\bf P}M_d^{simp}(Q)$ is
$$H_{T_Q}^{\ast}({\bf P}M_d^{simp}(Q))\simeq H_{T_Q}^{\ast}(T_Q/T)\simeq H_{T}^\ast(pt),$$
considered as a module over $H_{T_Q}^\ast(pt)$ via the restriction map. Thus, the equivariant cohomology $H_{T_Q}^{\ast}
({\bf P}M_d^{simp}(Q))$ is a torsion module
over $H_{T_Q}^\ast(pt)$, whereas the non-equivariant cohomology $H^{\ast}({\bf P}M_d^{simp}(Q)$ is two-dimensional, with
one-dimensional components in degrees $0$ and $1$.\\[2ex]
The inductive proof of Theorem \ref{euler} in section \ref{ind}
ultimately computes the Euler characteristic of ${\bf P}M_d^{simp}(Q)$ by
reduction to the case of single cycle quivers, thus detecting a special kind of grading on certain representations.
It is therefore natural to consider the class of simple
representation, which one might call string representations, defined as follows:

\begin{definition}\label{string} Suppose a primitive cycle
$$\omega\, :\, i_0\stackrel{\alpha_1}{\rightarrow}i_1\stackrel{\alpha_2}{\rightarrow}\ldots\stackrel{\alpha_s}{\rightarrow}
i_s=i_0$$ of dimension vector $d$ in $Q$ is given.
Define a simple representation $S_\omega$ of $Q$ of dimension vector $d$ as follows:\\[1ex]
For each vertex $i\in Q_0$, define
$${\cal K}_i=\{k=0,\ldots,s-1\, :\, i_k=i\}.$$
Consider a $d_i$-dimensional vector space $V_i$ with basis elements $b_k$ for $k\in {\cal K}_i$.
For each arrow $(\alpha:i\rightarrow j)\in Q_1$ and each $k\in {\cal K}_i$, define
$$(S_\omega)_\alpha(b_k)=\left\{\begin{array}{ccc}b_{k+1}&,&\alpha=\alpha_{k+1},\\ 0&,&\mbox{otherwise}.
\end{array}\right\}$$
\end{definition}
The reduction procedure of section \ref{ind} morally
means that the finitely many representations $S_\omega$ for primitive cycles of weight $d$ are the only ones
``responsible'' for (some aspects of) the global topology of ${\bf P}M_d^{simp}(Q)$.


\begin{thebibliography}{99}

\bibitem{A} M.~Artin: {\sl On Azumaya algebras and finite dimensional representations of rings.}
J. Algebra 11 (1969) 532--563.
 
\bibitem{BB} A.~Bia\l ynicki-Birula: {\sl Some theorems on actions of algebraic groups.} 
Ann. of Math. (2) 98 (1973), 480--497.

\bibitem{BS} R.~Bocklandt, S.~Symens: {\sl The local structure of graded representations.} Preprint 2005.

\bibitem{CG} N.~Chriss, V.~Ginzburg: {\sl Representation theory and complex geometry.}
Birkh\"auser Boston, Inc., Boston, MA, 1997.

\bibitem{EM} S.~Evens, I.~Mirkovi\'c: {\sl Fourier transform and the Iwahori-Matsumoto involution.} 
Duke Math. J. 86 (1997), no. 3, 435--464.

\bibitem{Gi} V.~Ginzburg: {\sl Principal nilpotent pairs in a semisimple Lie algebra. I.}
Invent. Math. 140 (2000), no. 3, 511--561.

\bibitem{GKM} M.~Goresky, R.~Kottwitz, R.~Macpherson: {\sl Equivariant cohomology, Koszul duality,
and the localization theorem.} Invent. Math. 131 (1998), no. 1, 25--83.  

\bibitem{Kl} A.~A.~Klyachko: {\sl Stable bundles, representation theory and Hermitian operators.}
Selecta Math. (N.S.) 4 (1998), no. 3, 419--445.

\bibitem{LBN} L.~Le Bruyn: {\sl Optimal filtrations on representations of finite-dimensional algebras.}
Trans. Amer. Math. Soc. 353 (2001), no. 1, 411--426.

\bibitem{LBP1} L.~Le Bruyn, C.~Procesi: {\sl Semisimple representations of quivers.} Trans. Amer. Math. Soc. 317 (1990),
no. 2, 585--598.

\bibitem{M} D.~Mumford, J.~Fogarty, F.~Kirwan: {\sl Geometric invariant theory. }
Third edition. Ergebnisse der Mathematik und ihrer Grenzgebiete (2), 34. 
Springer-Verlag, Berlin, 1994.

\bibitem{RCRP} M.~Reineke: {\sl Counting rational points of quiver moduli.} Preprint 2005. math.AG/0505389 

\bibitem{Reut} C.~Reutenauer: {\sl Free Lie algebras.}
London Mathematical Society Monographs. New Series, 7. Oxford Science Publications.
The Clarendon Press, Oxford University Press, New York, 1993.  


\end{thebibliography}
\end{document}